\documentclass[12pt,reqno]{amsart}
\usepackage{amssymb}
\usepackage{graphicx}
\usepackage{amsmath}
\usepackage{a4wide}
\usepackage{epsfig}

\newtheorem{theorem}{Theorem}

\newtheorem{corollary}[theorem]{Corollary}
\newtheorem{lemma}[theorem]{Lemma}
\newtheorem{proposition}[theorem]{Proposition}

\theoremstyle{remark}
\newtheorem{remark}[theorem]{Remark}

\def\leq{\leqslant}
\def\geq{\geqslant}
\def\al{\alpha}
\def\be{\beta}
\def\ga{\gamma}
\def\la{\lambda}

\def\tC{\mathbb C}
\def\G{\mathbb G}
\def\Gmn{{\mathbb G}_m^n}

\def\R{\mathbb R}
\def\Q{\mathbb Q}
\def\Z{\mathbb Z}

\def\eps{\varepsilon}

\def\om{\omega}

\def\th{\theta}
\def\bpsi{\boldsymbol\psi}
\def\cc {,\dots,}
\def\i{\mathbf i}
\def\bi{\i }
\def\j{\mathbf j}
\def\bom{\boldsymbol \omega}
\def\bla{\boldsymbol \lambda}
\def\bt{\mathbf t}
\def\bx{\mathbf x}

\def\ba{\mathbf a}
\def\tC{\tt C}
\def\C{\mathbb C}
\def\l{ \ell}
\def\bl{\boldsymbol \ell}
\def\bth{\boldsymbol \theta}
\def\brho{\boldsymbol \rho}
\def\lcm{\operatorname{lcm}}
\def\PN{P^{\phantom{\! N}}_{\! N}}

\begin{document}
\title[Salem numbers]{There are Salem numbers of every trace}
\author{ James McKee}
\address{Department of Mathematics, Royal Holloway\\University of London, Egham
Hill, Egham\\ Surrey TW20 0EX, UK}
\email{James.McKee@rhul.ac.uk}
\author{ Chris Smyth}
\address{School of Mathematics \\
University of Edinburgh\\
James Clerk Maxwell Building\\
King's Buildings, Mayfield Road\\
Edinburgh EH9 3JZ\\
Scotland, U.K.}
\email{C.Smyth@ed.ac.uk}
\date{19 Feb 2004}

\begin{abstract}
We show that there are Salem numbers of every trace. The nontrivial part of this result is for
 Salem numbers
of negative trace. The proof has two main ingredients. The first is a novel construction,
using pairs of polynomials whose zeros interlace on the unit
circle,
of polynomials of specified negative trace having one factor a Salem polynomial, with any other factors being cyclotomic.
 The second is an upper bound for
the exponent of a maximal torsion coset of an algebraic torus in a variety defined over the rationals.
 This
second result, which may be of independent interest, enables us to refine our construction to
avoid getting cyclotomic factors, giving
 a Salem polynomial of any specified trace, with a trace-dependent bound for
its degree.

We show also how our interlacing construction can be easily adapted to produce Pisot polynomials,
giving a simpler, and more explicit, construction for Pisot numbers of arbitrary trace than previously
known.
\end{abstract}

\maketitle

\section{\bf{Introduction}}

A {\it Salem number} is an algebraic integer greater than $1$ whose other
conjugates all lie in the closed disc $|z|\leqslant 1$, with at least one on $|z|=1$. Our
main result is the following.

\begin{theorem}\label{ThS} For every negative integer $-T$ there is a Salem
number
of trace $-T$ and degree at most $\exp\exp{(22+4T\log T)}$ .
\end{theorem}

It is easy to produce Salem numbers of any
nonnegative trace,   so the title of the paper
 is justified. The interest in this result is that, until now, all
 Salem numbers found had trace no smaller than $-1$ (\cite{Sm2}).  Furthermore, it is now known that
 a Salem number of degree $d\geq 10$ has trace at least
 $\lfloor1- d/9\rfloor$, and
 it seemed conceivable that there was a finite lower bound for the
 trace. For more details see the end of the paper.

To provide a little background, we  give a brief sketch of some facts about
Salem numbers. The minimal polynomial $P(z)$ of a Salem number $\tau$ is {\it reciprocal},
that is, it satisfies $z^{\deg P}P(1/z)=P(z)$, so that $\tau^{-1}$ is a
conjugate of $\tau$, and the coefficients of $P$ are ``palindromic''. All
conjugates of $\tau$ apart from $\tau$ and $\tau^{-1}$ lie on $|z|=1$,
and $P(z)$ has even degree. For
every $\eps>0$ and Salem number $\tau$ there is a $\la\in\Q(\tau)$ such that for $k=0,1,2\dots$
all $\la\tau^k$ have distance at most $\eps$ from an integer.
If a number field $K$ contains a Salem number $\tau$ of full degree $[K:\Q]$
then every full degree Salem number in $K$ is a power of the smallest such Salem
number in $K$. It is not known  whether there are Salem numbers arbitrarily
close to $1$. If ``Lehmer's conjecture" is true, then there are not. The
smallest known Salem number $1.176280818\cdots$, discovered by Lehmer in 1933, has minimal
polynomial $L(z)=z^{10}+z^9-z^7-z^6-z^5-z^4-z^3+z+1$. The polynomial $L(-z)$
had just appeared (in 1932) in Reidemeister's book \cite{R} as the Alexander
polynomial of a pretzel knot. For recent connections with knot
theory, see E. Hironaka \cite{H}.
The polynomial $L(z)$ can also be obtained
from the characteristic polynomial $E_{10}(x)$ of the (adjacency matrix of the) Coxeter graph $E_{10}=$
\leavevmode
\hbox{%
\epsfxsize0.8in \epsffile{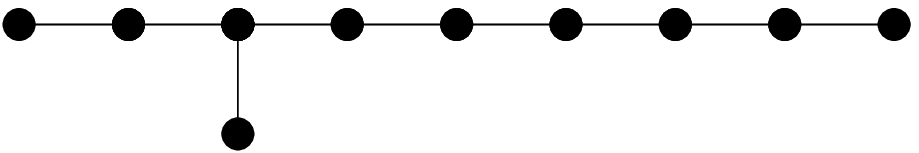}}
 by
the transformation $L(z)=z^5E_{10}(z^{1/2}+z^{-1/2})$. There are currently
$47$ known Salem numbers less than $1.3$
 (see Mossinghoff \cite{M}).

Salem numbers are closely related to Pisot numbers, which are much better
understood. A {\it Pisot number} is an algebraic integer
greater than $1$ whose other conjugates all lie in the open disc $|z|<1$.
For every Pisot number $\th$, the distance of $\th^n$ from the nearest integer tends to $0$ as $n\to\infty$.
The set of all Pisot numbers is a closed subset of the real line. Every
Pisot number is a limit point of Salem numbers, and
Boyd \cite[p.~327]{Bo1} has conjectured that Pisot numbers are the only limit points of Salem numbers.
If this is true, then the set of all Pisot and Salem numbers is also closed.
See Bertin {\it et al} \cite{BDGPS},
Boyd \cite{Bo1,Bo2}, Ghate and Hironaka \cite{GH} and Salem \cite{Sa} for these and other results about Salem and Pisot numbers.

It is already known (see \cite{MRS,McK}) that there are Pisot numbers of every trace.
However, we can greatly reduce the known upper bound for the smallest degree of a
Pisot number of given negative trace.

\begin{theorem}\label{ThP} For every negative integer $-T$ there is a Pisot number
of trace $-T$ and degree at most the sum of the first $2T+4$ primes.
\end{theorem}

This sum is asymptotic to $2T^2\log T$.
 The simple examples $z^3-z-1$, $z^2-z-1$ and $z-n\,(n\geq 2)$  of minimal polynomials of
Pisot numbers then show that there are Pisot numbers of every trace.

Computations for  negative trace down to $-25$ (see Section \ref{polycomp}) indicate that the upper bound on the degree
in Theorem \ref{ThS} should be comparable with that in Theorem \ref{ThP}. However, a proof of this
does not seem within reach at present.

In \cite{MRS}, for infinitely many degrees $d$  the existence of a Pisot number of degree $d$ and
 trace $<\frac{-\log d}{4(\log\log d)^{3/2}}$ was proved. Theorem \ref{ThP} improves this bound
 to $-c\sqrt{d}/\log d$ for some positive constant $c$.

One ingredient needed for the proof of Theorem \ref{ThS} is a result concerning
 the exponent of maximal torsion cosets on a variety (Theorem \ref{Th-Mbound}), which may be of independent
 interest.

To end the introduction, we mention one immediate consequence of Theorem \ref{ThS}.

\begin{corollary} For infinitely many  $n$ there is a totally
positive algebraic integer of degree $n$ and trace less than
$2n-\frac14\log\log n/\log\log\log n$.
\end{corollary}

This follows easily from Theorem \ref{ThS}, using the fact that
$\tau+1/\tau+2$ is totally positive for any Salem number $\tau$.

\section{\bf{Outline of the proof}}
There are two main ingredients in the proof of Theorem \ref{ThS}. The first is a new construction for
Salem numbers, which uses pairs of
polynomials whose zeros interlace on the unit circle. It is an extension of
the Salem number construction method used in \cite{MRS}, where the interlacing
polynomials arose from star-like trees. This new construction produces  a polynomial of
any specified negative trace that is, up to a possible cyclotomic factor, the
minimal polynomial of a Salem number (or a reciprocal Pisot number).

The purpose of the second ingredient is to get rid of the  possibility of a cyclotomic factor, while
at the same time bounding the degree of the Salem number.
It is based on ideas of Schmidt \cite{Sc}, and gives an upper
bound for the exponent of a maximal torsion coset on a variety.
This result is applied to a particular hypersurface to prove that the parameters in our
polynomial construction can be chosen so that the polynomial in fact has no
cyclotomic factor. This gives us our Salem number of the specified trace, with a bound on its degree.

\section{\bf{Construction of Salem and Pisot numbers by interlacing}}\label{inter}

\begin{lemma}\label{interlace-L} Suppose that
$\gamma>0$,
$\al_1<\be_1<\al_2<\cdots<\be_{d-1}<\al_d\leq A$, and
\begin{equation}\label{Ehx}
f(x):=\frac{\gamma\prod_{j=1}^{d-1}(x-\be_j)}{\prod_{j=1}^{d}(x-\al_j)}.
\end{equation}
Then $f(x)$ can be written as
\begin{equation}\label{Epos-lam}
f(x)=\sum_j\frac{\la_j}{x-\al_j},\quad\text{with } \la_j>0\quad\text{for all } j. 
\end{equation}

Further, the equation $f(x)=1$ has real roots $\ga_1,\dots,\ga_d$, where
$\al_1<\ga_1<\be_1<\al_2<\ga_2<\be_2\cdots<\ga_{d-1}<\be_{d-1}<\al_d<\ga_d$. Also $\ga_d>A$ if
and only if $f(A)>1$.

Conversely, every $f(x)$ of the form (\ref{Epos-lam}) can be written in the form
(\ref{Ehx}) for some $\gamma>0$ and $\be_1\cc  \be_{d-1}$ that interlace with the $\al_j$.
\end{lemma}

\begin{proof} The interlacing condition for the roots easily implies
(\ref{Epos-lam}). Then the results follow immediately on applying the Intermediate Value Theorem
to $\sum_j\frac{\la_j}{x-\al_j}$.
\end{proof}

We say that a pair of relatively prime polynomials $p$ and $q$ satisfy the {\em circular interlacing condition} if
they both have real coefficients, positive leading term,  and all their zeros lie on the unit circle,
 and interlace there. This last condition means that,
 progressing clockwise around the unit circle,  a zero of $p$ and a zero of $q$
are encountered alternately. Thus $p$ and $q$ have the same degree, and neither has a multiple zero. Note too that if  $p$ and $q$ satisfy the  circular interlacing condition, so do
$p(z^n)$ and $q(z^n)$ for $n=1,2,\dots$. In particular, the pair $z^n-1$ and $z^n+1$ satisfy it.

  Pisot numbers whose minimal polynomials are reciprocal behave in some ways like Salem numbers.
It is clear that they must be quadratic.

\begin{proposition}\label{interlace-P}  Suppose that the
polynomials $p$ and $q$ satisfy the  circular interlacing condition, have integer coefficients,
and that $p$ is monic (and thus cyclotomic). Then

\begin{enumerate}

\item[(a)] if $p(1)=0$, or $q(1)=0$ and $2p(1)-q'(1)<0$, then $(z^2-1)p(z)-zq(z)$ is the minimal polynomial
of a Salem number (or perhaps a reciprocal Pisot number), possibly multiplied by a cyclotomic polynomial.
{\it [Note: one of $p(1)$ and $q(1)$ is always $0$.]}

\item[(b)]  always $(z^2-z-1)p(z)-zq(z)$ is the minimal polynomial of a Pisot number.

\end{enumerate}

\end{proposition}

\begin{proof} Firstly, it is clear that, as the zeros of $p$ and $q$ interlace, both $1$ and
$-1$ must be zeros of  $pq$,
all other zeros of both $p$ and $q$ occurring in complex conjugate pairs. Put
$z+1/z=x$, real and in $[-2,2]$ for $z$ on
the unit circle.

\begin{enumerate}
\item[(a)]
Let $\gamma$ be the leading coefficient of $q$. 

Suppose first that $p$ and $q$ have even degree $2d$. If $z^2-1$ divides $q$,
then
\begin{equation}\label{Efdef}
f(x):=\frac{z}{z^2-1}\cdot
\frac{q(z)}{p(z)}=\frac{\gamma\prod_{j=1}^{d-1}(x-\be_j)}{\prod_{j=1}^{d}(x-\al_j)},
\end{equation}
where the $\al_j=r_j+1/r_j$, $\be_j=s_j+1/s_j$ for zeros $r_j$ of $p$, $s_j$ of
$q$, and
$$
-2<\al_1<\be_1<\al_2<\cdots<\be_{d-1}<\al_d< 2.
$$

Thus, by Lemma \ref{interlace-L},
 this quotient is equal to $\sum_j\frac{\la_j}{x-\al_j}$ for some $\la_j>0$. 
 On
the other hand, if $z^2-1$ divides $p$, then
\begin{equation*}
\frac{z}{z^2-1}\cdot
\frac{q(z)}{p(z)}=\frac{\gamma\prod_{j=1}^{d}(x-\be_j)}{(x^2-4)\prod_{j=1}^{d-1}(x-\al_j)},
\end{equation*}
with
$$
-2<\be_1<\al_1<\be_2<\cdots<\al_{d-1}<\be_d< 2.
$$
Thus Lemma \ref{interlace-L} can be applied again to give the same conclusion.

Now suppose that $p$ and $q$ have odd degree $2d+1$. Then for $\eps$ equal to one of $\pm 1$,
$(z-\eps)$ divides $q$ and $(z+\eps)$ divides $p$.
Then
\begin{equation*}
\frac{z}{z^2-1}\cdot
\frac{q(z)}{p(z)}=\frac{\gamma\prod_{j=1}^{d}(x-\be_j)}{(x+2\eps)\prod_{j=1}^{d}(x-\al_j)},
\end{equation*}
and
 $$
-2<\be_1<\al_1<\be_2<\cdots<\be_d<\al_d< 2 \quad\text{for $\eps=1$}
$$
$$
-2<\al_1<\be_1<\al_2<\cdots<\al_d<\be_d <2\quad\text{for $\eps=-1$}
$$
so that again Lemma \ref{interlace-L} applies.

Now one of $\{p(1),q(1)\}$ is zero and the other positive, and $(z^2-1)p(z)-zq(z)$ is the numerator
of $1-f(x)$. Hence the conditions given at $z=1$ are clearly
those necessary and sufficient for $(z^2-1)p(z)-zq(z)$
to have a real zero greater than $1$.

\item[(b)]

We consider the sum 
\begin{equation*}
\frac{q(z)}{p(z)}+\frac{z^n+1}{z^n-1}
\end{equation*}
and write its uncancelled numerator and denominator as
$q^*(z)g(z)$ and $p^*(z)g(z)$ respectively, where $q^*$ and $p^*$
are relatively prime. Then, because the pairs $\{q,p\}$ and
$\{z^n+1,z^n-1\}$ satisfy the circular interlacing condition, so
do $q^*$ and $p^*$, by Proposition \ref{Psum}. (Note that there is
no circularity, as the proof of Proposition \ref{Psum} uses only
part
 (a) of this proposition.) Then
$g(z)((z^2-1)p^*-zq^*)=z^n((z^2-z-1)p-zq)-(z^2+z-1)p+zq$ is of the
form $z^nR(z)\pm R^*(z)$, where $R(z)=(z^2-z-1)p-zq$ and
$R^*(z)=z^{\deg R} R(1/z)$, using the fact that one of $p^*$ and
$q^*$ (say $p^*$) is reciprocal, and the other (say $q^*$)
satisfies $z^{\deg q^*} q^*(1/z)=-q^*(z)$.
 Now, for any $\eps>0$, apply Rouch\' e's Theorem on
the circle $|z|=1+\eps$, and let $n\to\infty$. This shows that
$R(z)=(z^2-z-1)p-zq$ has at most one zero in $|z|>1$. Also, $R$
has no zeros on $|z|=1$, as any such zero would also be a zero of
$R^*$, and, as is easily checked, $R$ and $R^*$ are relatively
prime. Finally, because one of $\{p(1),q(1)\}$ is zero and the
other positive, in fact $R(1)<0$, so that $R$ does have one zero
on $z>1$.

\end{enumerate}

\end{proof}

\begin{proposition}\label{Psum}
Suppose that the pairs of polynomials $p_i,q_i\quad(i=1\cc I)$ each satisfy the circular interlacing
condition. Then  $\sum_i \frac{q_i(z)}{p_i(z)}$
is equal to a quotient $\frac{q(z)}{p(z)}$, where  $p$
and $q$ also satisfy the circular interlacing
condition. 
\end{proposition}

\begin{proof} From Lemma \ref{interlace-L} and the  proof of Proposition \ref{interlace-P}(a), we know that for each $i$
\begin{equation}
\frac{z}{z^2-1}\cdot
\frac{q_i(z)}{p_i(z)}=\sum_j\frac{\la_j}{x-\al_j}
\end{equation}
where $x=z+1/z$, the $\la_j$ are positive, and the $\al_j$ are all
real and in $[-2,2]$. On adding, the same applies to
$\sum_i\frac{z}{z^2-1}\cdot\frac{q_i(z)}{p_i(z)}$.
Hence, by Lemma
 \ref{interlace-L}, this sum is equal to a positive scalar multiple of a quotient of polynomials
$\frac{\prod_{j=1}^{d-1}(x-\be_j)}{\prod_{j=1}^{d}(x-\al_j)}$, where
$-2\leq\al_1<\be_1<\al_2<\cdots<\be_{d-1}<\al_d\leq 2$. Then on substituting $x=z+1/z$ and considering
separately the cases when  $\al_1=-2$ or $\al_d= 2$, we get the main result.

\end{proof}

\section{\bf{The exponent of maximal torsion cosets}}
As usual, let $\G_m$ denote the multiplicative group of $\C$. An
{\it $r$-dimensional   subtorus} $H$ of $\G_m^{n}$ is a subgroup
of the group $\G_m^{n}=\{(x_1\cc  x_n)\mid x_i\ne 0\}$ where, for
some $r$, parameters $t_1\cc  t_r$ and integer matrix
$E=(e_{ji})_{ (j=1\cc  r;i=1\cc  n)}$ of rank $r$ we have
$x_i=t_1^{e_{1i}}\cdots t_r^{e_{ri}}$. It is an algebraic subgroup
of $\G_m^{n}$, defined by the equations $\{\bx^\ba=1\mid \ba\in
A\}$, where the $\ba\in A$ span the lattice of integer vectors
orthogonal to the rows of $E$. A {\it torsion coset} is a
translate $\boldsymbol \om H$ of $H$ by a torsion point
$\boldsymbol  \om=(\om_1\cc  \om_n)$, the $\om_i$ being roots of
unity. An {\em exponent} of $\boldsymbol \om H$ is any multiple of
its order as an element of the group $\G_m^{n}/H$. A {\it maximal
torsion coset} of a variety $V$ is a torsion coset not properly
contained in any other torsion coset in $V$. Results of Laurent
\cite[Th. 2]{L}, Bombieri and Zannier \cite[Th. 2]{BZ}, and
Schmidt \cite[pp. 159--60]{Sc} state that for any variety
 $V\subset \G_m^{n}$ defined over a number field $K$,
the union of all torsion cosets contained in $V$ is in fact contained in a union
of finitely many maximal torsion cosets in $V$, with an upper bound for this number depending only on
the parameters of $K$ and $V$. Furthermore, in \cite{Sc} Schmidt has
given an explicit bound of this kind.

The finiteness of the number of maximal torsion cosets in $V$ immediately implies
the existence of a single exponent for all these cosets. This fact
can be used to prove, as in Section \ref{S-Proof-Th1}, that there are Salem numbers of a given trace, but
without the upper bound on the smallest degree of such a number. The results that
follow (Corollary \ref{C-expo} in particular) are needed to produce this degree bound.

We denote a typical torsion coset by
${\tC}=\boldsymbol \om \mathbf t^E=(\om_i\prod_{j=1}^r{ t_j}^{e_{ji}})_{(i=1\cc  n)}\subset\Gmn$,
 $E$ being an $r\times n$ integer matrix of rank $r$.

 Consider a system of linear
equations
\begin{equation}\label{Eu}
\sum_{i=1}^N a_{\ell i}u_i=0\quad(\ell=1,\dots, L ).
\end{equation}
Following Schmidt, a solution $\mathbf u=(u_1,\dots,u_{\! N})\in\G_m^N$ will be called {\it nondegenerate} if there is
no subset $I$ of $\{1,\dots,N\}$ with
$0<\# I<N$ such that
$$
\sum_{i\in I} a_{\ell i}u_i=0\quad(\ell=1,\dots, L ).
$$

\begin{lemma}\label{rootorder}(see \cite[p. 168--9]{Sc}, \cite{CJ} ) Suppose we have a nondegenerate
solution of
(\ref{Eu}) where the $u_i$ are all roots of unity. Then, up to a factor of
proportionality, the $u_i$ are all $\PN$-th roots of unity,
where $\PN$ is the product of all primes up to $N$.
\end{lemma}

In fact, their result tells us that such solutions are $m$-th roots of unity, where $m$ is the product
of at most $2\sqrt{N}$ distinct primes
$p\leq N$. However, we need an exponent valid uniformly for solution sets of  different such
$N$-term equations. This is
why we take  $\PN$-th roots of unity, $\PN$ being the $\lcm$ of all such $m$. A uniform `killer' exponent
is provided by the following result, and its corollary.

\begin{theorem}\label{Th-Mbound} Suppose that $V$ is an affine variety in $\G_m^n$ defined over $\Q$,
given say by polynomial
equations
\begin{equation}\label{EV}
\sum_\bi a_{\ell \i } \bx^{\i }=0 \qquad(\ell=1,\dots, L )
\end{equation}
with total degree $d$. Suppose also that
the set
$$
\mathcal N(V)=\{\i \in \Z^n \mid a_{\ell \i } \ne 0 \text{ for some }\ell\}
$$
has diameter $D(V)$.
Then every
$(n-k)$-dimensional maximal torsion  coset on $V$ has an exponent $m\PN$ for some integer
$m\leq D(V)^{2k}k^{k/2}$.
 Here $N:=\# \mathcal N(V)\leq \binom{n+d}{d}$.
\end{theorem}

\begin{proof} The ingredients for the proof come from Schmidt \cite{Sc}.
Take $r=n-k$ and a maximal $r$-dimensional torsion coset ${\tC}=\boldsymbol \om \mathbf t^E$ on $V$,
so that
$$
\sum a_{\ell \bi}\boldsymbol \om^\i \mathbf t^{E\i }=0\quad(\ell=1,\dots, L ).
$$
Our aim is to find $\bom_1$ with also ${\tC}=\bom_1 \mathbf t^E$, with $\bom_1$ a vector of
 $(m\PN)$-th roots of unity for some $m\leq D(V)^{2k}k^{k/2}$.
Now for any $\mathbf j\in \Z^r$ the coefficient of $\bt^\j$ is 

\begin{equation}\label{om-eqn}
\sum_{\i :E\i =\mathbf j }a_{\ell \bi}\boldsymbol \om^\mathbf
i=0\quad(\ell=1,\dots, L ).
\end{equation}
Here the sums over $\i $ are taken over
all relevant $\i $ in
$\mathcal N(V)$.
Now (\ref{om-eqn}) may
be degenerate, splitting into nondegenerate equations

\begin{equation}\label{om-eqn-k}
\sum_{\i \in I_q}a_{\ell \bi}\boldsymbol \om^\i =0\quad(\ell=1,\dots, L ,q\in
Q \hbox{ say})
\end{equation}
for nonempty subsets $I_q$ of $\Z^n$.
Now, for a single $q$,  apply Lemma \ref{rootorder} to (\ref{om-eqn-k}), to obtain

\begin{equation}\label{om-eqn-q}
\sum_{\i \in I_q}a_{\ell \bi}\boldsymbol \om^{\i -\mathbf
i_q}=0\quad(\ell=1,\dots, L )
\end{equation}
where $\i _q$ is some fixed vector in $I_q$. Here, the number of terms is at most $N$  .
Then for all $q\in Q$, we have from Lemma \ref{rootorder} that  all $\boldsymbol \om^{\i -\bi_q} (q\in
Q) $ are vectors of
 $\PN$-th roots of unity.

Recalling that $k=n-r$, we claim that the set of all vectors $\{\i -\i _q\mid \i \in
I_q,q\in Q\}$ generates a $k$-dimensional sublattice ${\mathcal L}_{\tC}$ of $\Z^n$.
 For, from (\ref{om-eqn}), the lattice ${\mathcal L}^E$ spanned by the rows of $E$ is orthogonal
to ${\mathcal L}_{\tC}$, and so ${\mathcal L}_{\tC}$
has dimension $\leq n-r$. But if the inequality were strict,
there would be a nonzero vector $\i '\in\Z^n$ orthogonal  to ${\mathcal L}_{\tC}$ and
not in the rational span of ${\mathcal L}^E$.   Then
for $\i \in I_q$ we would have $\i ' \cdot \i
= \i ' \cdot \i _q$
(and also of course $E{\i }=E{\i _q}$), so for any $u\in\G_m$ we would have,
for $\ell=1,\ldots,L$,
 \begin{align*}
     \sum a_{\l\i } {\boldsymbol \om}^{\i }
\mathbf t^{E\i } u^{\i '\cdot\i }&=
\sum_q\sum_{\i \in I_q} a_{\l\i }{\boldsymbol \om}^{\i } \mathbf t^{E\i _q} u^{\i '\cdot \i _q}    \\
    &=
\sum_q \mathbf t^{E\i _q}u^{\i '\cdot \i _q}
\sum_{\i \in I_q} a_{\l\i }{\boldsymbol \om}^{\i }
\\
    &=0 \,,
\end{align*}
 and so the larger torsion coset $\boldsymbol \om \mathbf {t'}^{E'}$ would lie on $V$, where $\mathbf
t'=(\mathbf t,u)$ and $E'=\binom{E}{\i '}$, contradicting the maximality of
 $\boldsymbol \om \mathbf t^E$.

Next take a basis  $\bl_1\cc  \bl_k$ of vectors in $\{\i -\i _q\mid \i \in
I_q,q\in Q\}$ for $\mathcal L_{\tC}$, and put
$\bom=e^{i\bth}=(e^{i\th_1}\cc e^{i\th_n})$.
 Write
$\bth=\sum_{j=1}^k \la_j\bl_j+\bpsi$, where  $\bpsi=\brho E$ for some $\brho\in\R^r$.
Then, on solving the  system of linear equations
$$
\bl_i\cdot \bth=\sum_{j=1}^k \la_j\bl_i\cdot\bl_j\quad (i=1\cc  k)
$$
and using the fact that $\PN\bl_i\cdot \bth\equiv 0\pmod {2\pi}\quad(i=1\cc k)$, we see that
$\PN\det(\bl_i\cdot\bl_j)\bla\equiv 0\pmod {2\pi}$. Note too that $\det(\bl_i\cdot\bl_j)\ne 0$.
Then, using the Cauchy-Schwartz and
Hadamard inequalities, we have that
$$
|\det(\bl_i\cdot\bl_j)|\leq D(V)^{2k}k^{k/2}.
$$
Put
 $\bt_1=e^{-i\brho}$. Then $\bt_1^E=e^{-i\bpsi}$, and for $\bom_1=\bom\bt_1^E$ and some
 $m'=m \PN$ with $m\leq{D(V)}^{2k}k^{k/2}$
 we have $\bom_1^{m'}=1$.
Since
$$
{\tC}=\bom{\bt}^E=\bom({\bt}_1{\bt})^E=\bom_1{\bt_2}^E
$$
say, we see that $\tC$ has exponent $m'$.
\end{proof}

This result immediately gives us a killer exponent $K$ valid for all maximal torsion cosets on $V$.

\begin{corollary}\label{C-expo} Let $V$ be as in the Theorem, and $K=\PN\lcm(1,2,\dots, D(V)^{2n}n^{n/2})$,
 where $N=\# \mathcal N(V)$. Then every maximal torsion coset of $V$ has exponent $K$.
\end{corollary}

\section{Maximal torsion cosets on a particular hypersurface}

We shall be applying the results of the previous section to the affine hypersurface $h(\bx)=0$, where
$\bx=(x_0\cc x_{n})\in\G_m^{n+1}$ and
$$
h(\bx)=2(x_0^2-1)\prod_{i=1}^n(x_i-1)-x_0\sum_{j=1}^n(x_j+1)\prod_{\underset{i\ne j}{i=1}}^n(x_i-1).
$$
The reason for looking at this hypersurface is that we shall apply the identity
\begin{equation}\label{h-inter}
 \frac{h(\bx)}{2 x_0\prod _{i=1}^n(x_i-1)}=
 \frac{x_0^2-1}{x_0}-\frac {1}{2}\sum_{i=1}^n\frac{x_i+1}{x_i-1},
 \end{equation}
  which is connected to our interlacing considerations of Section \ref{inter}.
\begin{lemma} The only maximal torsion cosets of $h$ with
$x_0$ nonconstant are the algebraic subgroups $B_{ij}$ of $\G_m^{n+1}$, where $i\ne j$ are both
nonzero, and
$$
B_{ij}=\{\mathbf  x\mid x_i=x_j=1,x_0=t_0,x_\ell=t_\ell(\ell\ne i,j)\},
$$
of rank $n-1$.
\end{lemma}
\begin{proof} Clearly no point on $h=0$ can have just one $x_i=1.$ If $\mathbf  x$
with any two
$x_i,x_j=1$ is on $h=0$ then it belongs to $B_{ij}$. Thus any other rank $r$
maximal
torsion coset with $x_0$ nonconstant has no $x_i$ identically $1$, so that we must have
$x_0=\om_0 t_1^{e_{10}}\cdots t_r^{e_{r0}}$, where $(e_{10}\cc  e_{r0})\ne 0$ and
$x_i=\om_i t_1^{e_{1i}}\cdots t_r^{e_{ri}}$, where $(e_{1i}\cc  e_{ri})\ne 0$
 whenever $\om_i=1$. By avoiding certain hyperplanes we can choose $\pm (k_1\cc k_r)
 \in \Z^r$ not orthogonal to any nonzero $(e_{1i}\cc  e_{ri})$. Then
 for $(t_1\cc  t_r)=(t^{k_1}\cc  t^{k_r})$, $x_i=\om_it^{\ell_i}$ where
 $\ell_i:=\sum_{j=1}^r k_je_{ji}\ne 0$ when $\om_i = 1$, and, by choice of the sign,   $ \ell_0>0$. Now as
 $t\to\infty$,
   the right-hand side of  (\ref{h-inter}) goes to infinity, so that the coset cannot be on $h=0.$
 \end{proof}

 We now estimate the killer exponent $K$, valid for every maximal torsion coset on this hypersurface,
    defined over $\G_m^{n+1}$.

 \begin{lemma}\label{L-estimate} There is a killer exponent $K$ with $ \log\log K<0.2+(3(n+1)/2)\log(n+3)$
 for the hypersurface $h=0$. Further,
 $K$ can be chosen with all its prime factors less than $(n+3)^{3(n+1)/2}$.
\end{lemma}
\begin{proof}
      The hypersurface has diameter  $D=\sqrt{n+4}$, degree $d=n+2$,
 $N=\# \mathcal N(h=0)=3\cdot 2^n$.
  Hence $D^{2n+2}(n+1)^{(n+1)/2}< (n+3)^{3(n+1)/2}$, and Corollary \ref{C-expo} gives
  $K=\PN\!\cdot\lcm(1,2,\dots,D^{2n+2}(n+1)^{(n+1)/2}) $, with all prime factors of $K$
  less than $(n+3)^{3(n+1)/2}$.
  Then, using standard bounds of Rosser and Schoenfeld \cite{RS} for
 the arithmetical functions $\theta,\psi$ we
 obtain
 \begin{align*}
    \log K &<  \theta(3\cdot 2^n)+\psi((n+3)^{3(n+1)/2})  \\
    &< 1.02\cdot 3\cdot 2^n+1.04\cdot (n+3)^{3(n+1)/2} \\
    &< 1.2(n+3)^{3(n+1)/2},
\end{align*}
giving the upper bound claimed.
 \end{proof}

 \section{\bf{Proof of Theorem \ref{ThS}}}\label{S-Proof-Th1}

 The following lemma will complete the proof of Theorem \ref{ThS}.
 \begin{lemma} For given even $n$ there are positive integers $k_1\cc  k_n$
 such that
 $$
 h(t,t^{k_1}\cc  t^{k_n})=2(t-1)^{n-1}S(t),
 $$
 where $S(t)\in\Z[t]$ is monic irreducible and the minimal polynomial of a Salem
 number of trace $T:=1-n/2$. Further, $S$ has degree less than $\exp\exp(22+4T\log T)$.
 \end{lemma}

 \begin{proof} For all maximal torsion cosets of $h$ with $x_0$ constant
 (ie all except the $B_{ij}$) we can suppose that the constant $x_0$-values
 are all $K$-th roots of unity, where furthermore $K$ has been chosen
 minimally. Note that $K$ is certainly even, because  the  point
 $x_0=x_1=\dots=x_{\! N}=-1$ lies on $h=0$ and, as it is on no $B_{ij}$, must lie in one
 of the constant-$x_0$ maximal torsion cosets. Take
 $k_1=K$, and $k_2\cc  k_n$ as the smallest $n-1$ primes not dividing $K$. Then
 all $k_1,\dots,k_n$ are pairwise relatively prime. We now assert that for every
 root of unity $\om$ and ${\om}^{\mathbf  k}=(\om,\om^{k_1}\cc  \om^{k_n})$
 with
 $h({\om}^{\mathbf  k})=0$, we have $\om=1$.

 For ${\om}^{\mathbf  k}$ belongs to some maximal torsion
 coset. If ${\om}^{\mathbf  k}$ has at least two components $=1$, then
(by the extended euclidean algorithm) $\om=1$. Alternatively,  it
 belongs to no $B_{ij}$, and so to some maximal torsion coset with $x_0$
 constant, $x_0=\om$, and so $\om^{ K }=\om^{k_1}=1$. This is impossible, as we cannot have just one $x_i=1$,
 as noted above. This proves the assertion.

 It is easy to check that $(d/dt)^nh(t,t^{k_1}\cc  t^{k_n})$ evaluated at $t=1$,
 is nonzero.  Furthermore, $h(t,t^{k_1}\cc
 t^{k_n})\equiv \prod _{i=1}^n(x_i-1)(2(t^2-1)-tn)\equiv 0\pmod 2$, so that all
 coefficients of $h$ are even. This gives the stated factorization $2(t-1)^{n-1}
 S(t)$ of
 $h(t,t^{k_1}\cc  t^{k_n})$.  Also, as all $k_i\geq 2$,
 $$
 h(t,t^{k_1}\cc  t^{k_n})=2t^{2+\sum k_i}-nt^{1+\sum k_i}+\dots
 $$
 showing that $S(t)$ has trace $1-n/2$.

 Finally, to show that $S$ is the minimal polynomial of a Salem number,
 observe that we have shown that none of its zeros are roots of unity.
 Now since $t^k+1$ and $t^k-1$ satisfy the circular interlacing condition, so does the sum
 $\frac{1}{2}\sum_{i=1}^n\frac{t^{k_i}+1}{t^{k_i}-1}$, by Proposition \ref{Psum}, so we can
 write it as
$q(t)/p(t)$, where $p$ and $q$ satisfy the circular interlacing
condition. Furthermore, as $n$ is even, $q$ has integer coefficients, as does $p=\prod_i(t^{k_i}-1)/(t-1)^{n-1}$.
 Hence,  as $p(1)=0$, the numerator $S(t)=(t^2-1)p(t)-tq(t)$ of the right-hand side of
 (\ref{h-inter})  with $x_0=t,x_i=t^{k_i}$ is, by Proposition \ref{interlace-P}(a), the minimal polynomial
 of a Salem number  of trace $1-n/2$.

Now the degree of $S$ is $2+\sum_ik_i-(n-1)$, and from Lemma \ref{L-estimate} we can take $k_2,\dots,k_n$ to
be the smallest $n-1$ primes greater than $(n+3)^{3(n+1)/2}$. By Bertrand's Postulate (Chebyshev's Theorem),
this gives $\deg S<K+(n-1)2^{(n-1)}\cdot (n+3)^{3(n+1)/2}<2K$, and $\log\log\deg S<\log\log K+\log 2/\log
K<0.2+(3(n+1)/2)\log(n+3)+0.1$. For $n=2T+2$ one readily checks that this is less than $22+4T\log T$.
 \end{proof}

\begin{remark} There are many  maximal torsion
 cosets with $x_0$ constant, for instance
  $x_0=-1,x_1=x_2^{-1}=t_1\cc  x_{n-1}=x_n^{-1}=t_{n/2}.$
  Also one can for instance construct some for $x_0=1$ using the
   identity $3\cot \pi/3-\cot\pi/6=0$.
  \end{remark}

 \section{\bf{Proof of Theorem \ref{ThP}}}
 \begin{proof} The proof is much easier for Pisot numbers,
  as there are no possible cyclotomic factors to dispose of.
 We replace the fraction $(t^2-1)/t$ in (\ref{h-inter})
  by $(t^2-t-1)/t$, and then can simply choose the parameters $k_i$  to be the first $n$  primes. Thus,
  again putting $\frac {1}{2}\sum_{i=1}^n\frac{t^{k_i}+1}{t^{k_i}-1}=q(t)/p(t)$, the polynomial
   $(z^2-z-1)p(z)-zq(z)$ will be the minimal polynomial of a Pisot number of trace $2-n/2$.
 \end{proof}

 \section{{\bf Computing Salem and Pisot numbers of negative trace.}}\label{polycomp}
 Salem and Pisot numbers of negative trace can be produced using
 $h(t,t^{k_1}\cc  t^{k_n})$, as in the previous sections. Thus,
 for the Pisot numbers of trace $-T$, the first $2T+4$ primes are used
 for the $k_i$. For the Salem numbers of trace $-T$, the first $2T+2$ primes are used
 for the $k_i$. In particular, $k_1$ is taken to be simply $2$, instead of the
 very large killer exponent $K$ used in the proof above. Computation using Maple shows that this produces
  a polynomial free of cyclotomic
 factors for  $T\leq 25$, giving a Salem number of trace $-T$ and degree equal to the sum of the first
 $2T+2$ primes minus $2T-1$ (for instance degree $5540$ for trace $-25$).
 However, we do not know whether
 this always happens. It would of course be nice if this could be proved, as we would then obtain a degree bound in
 Theorem \ref{ThS} as good as that in Theorem
 \ref{ThP}.

 Here is some pseudocode that gives the minimal polynomials.
 For a Salem number of trace $-T$:

 \bigskip

 {\tt $r=1;S=(z^2-1)(z-1);Q=z;$

 for $j=1,\dots,T+1$ do

 $q=\text{nextprime}(r);r=\text{nextprime}(q)$;

 $S=\frac{z^q-1}{z-1}\cdot \frac{z^r-1}{z-1}\cdot S-\frac{z^{q+r}-1}{z-1}\cdot Q;$

 $Q=\frac{z^q-1}{z-1}\cdot \frac{z^r-1}{z-1}\cdot Q$;

 enddo

 if $\gcd(S(z),S(-z)S(z^2)S(-z^2))=1$ then print$(S)$;

 endif}

 \bigskip

 The $\gcd$ condition tests for cyclotomic factors, and is based on the fact that a root of unity $\om$ is
 conjugate
 to one of $-\om$, $\om^2$ or $-\om^2$. See \cite{BS} for further developments of this idea.
 For instance, for $T=2$ we obtain the (reciprocal)  Salem polynomial
 {\tiny\begin{multline*}
 S_{-2}(z)=z^{38} + 2\,z^{37} - 2\,z^{36} -
19\,z^{35}
 - 57\,z^{34} - 123\,z^{33} - 222\,z^{32} - 357\,z^{31} -
527\,z^{30} - 727\,z^{29} - 950\,z^{28} - 1190\,z^{27} \\
\mbox{} - 1440\,z^{26} - 1692\,z^{25} - 1936\,z^{24} - 2161\,z^{
23} - 2355\,z^{22} - 2506\,z^{21} - 2602\,z^{20}
 - 2635\,z^{19} - 2602\,z^{18} - \cdots +1
\end{multline*}
}

Note that the Salem polynomials $S$ produced by this method have $|S(-1)S(1)|$ large. This is easily
seen by putting $z=\pm 1$ in the pseudocode. An
interesting question is whether there are Salem numbers with arbitrary trace and $|S(-1)S(1)|=1$,
the so-called {\it unramified} Salem numbers (see Gross and McMullen \cite{GM}).

  For a Pisot number of trace $-T$:

  \bigskip

 {\tt $r=1;P=z^2-z-1;Q=z;$

 for $j=1,\dots,T+2$ do

 $q=\text{nextprime}(r);r=\text{nextprime}(q)$;

 $P=(z^q-1)(z^r-1)\cdot P-(z^{q+r}-1)\cdot Q;$

 $Q=(z^q-1)(z^r-1)\cdot Q$;

 enddo

 print$(P)$;}

 \bigskip

Finally, we justify the statements in the Introduction. We note that there are Salem numbers of every nonnegative trace: for $n>0$ the polynomial
 $z^4 - nz^3 - (2n+1)z^2 - nz + 1=z^2((z+1/z)^2-n(z+1/z)-(2n+3))$ is easily seen to be
 the minimal polynomial of a Salem number of trace $n$. This follows from
 the fact that $x^2-nx-(2n+3)$ has one zero in $(-2,2)$ and the other zero greater than $2$.
  Also, $z^6-z^4-2z^3-z^2+1=z^3((z+1/z)^3-4(z+1/z)-2)$ is the
 minimal polynomial of a Salem number of zero trace.

 The lower bound $\lfloor 1-d/9\rfloor$ for the trace of a degree $d\geq 10\,$ Salem number
follows from the fact that the trace of a totally positive algebraic
integer of degree $n\geq 5$ is greater than $16n/9$ (\cite{MS}) on noting that for a Salem number $\tau$ of degree $d$ and trace $-T$,
the number $\tau+1/\tau+2$ is totally positive of degree $d/2$ and trace $d-T$.
Thus (turning the inequality around) for $-T\leq -2$ every Salem number of
trace  $-T$ has degree at least
$$
2\left\lceil \frac{9T}{2}\right\rceil +2=
\begin{cases}
18k+2 \text{ for } -T=-2k\\
18k+10 \text{ for } -T=-(2k+1).
\end{cases}
$$
This inequality is sharp for $T=2$ (\cite{MS}).
 The only
Salem number of degree less than $10$ having negative trace is the one with minimal polynomial
$z^8+z^7-z^6-4z^5-5z^4-4z^3-z^2+z+1$ (\cite{Sm2}).

\bigskip

{\it Acknowledgment.}
  We are grateful to the London Mathematical Society for their financial support for this work through a
  Collaborative Small Grant. We also thank the referee for helpful
  comments.

\end{document}